\documentclass{jnmp}

\usepackage{amsmath}

\JNMPnumberwithin{equation}{section}

\newtheorem{theorem}{Theorem}[section]
\newtheorem{proposition}{Proposition}[section]
\newtheorem{lemma}{Lemma}[section]
\newtheorem{conjecture}{Conjecture}[section]

\theoremstyle{definition}
\newtheorem{definition}{Definition}[section]
\newtheorem{example}{Example}[section]

\begin{document}

\renewcommand{\evenhead}{B~Agrebaoui, M~Ben Ammar, N~Ben Fraj and V~Ovsienko}
\renewcommand{\oddhead}{Deformations of Modules of Differential Forms}

\thispagestyle{empty}

\FirstPageHead{10}{2}{2003}{\pageref{ovsienko-firstpage}--\pageref{ovsienko-lastpage}}{Letter}

\copyrightnote{2003}{B~Agrebaoui, M~Ben Ammar, N~Ben Fraj and V~Ovsienko}

\Name{Deformations of Modules of Differential Forms}
\label{ovsienko-firstpage}

\Author{B~AGREBAOUI~$^\dag$, M~BEN AMMAR~$^\dag$, N~BEN FRAJ~$^\dag$
and V~OVSIENKO~$^\ddag$}

\Address{$^\dag$~D\'epartement de Math\'ematiques, Facult\'e des Sciences de Sfax,
Route de Soukra,\\
~~3018 Sfax BP 802, Tunisie\\
~~E-mail: bagreba@fss.rnu.tn, Mabrouk.BenAmmar@fss.rnu.tn,\\
~~$\phantom{\mbox{E-mail:}}$~nizar.benfraj@issato.rnu.tn\\[10pt]
$^\ddag$~Institut Girard Desargues, Universit\'e Lyon-I, 21 Avenue Claude Bernard,\\
~~69622 Villeurbanne Cedex, France\\
~~E-mail: ovsienko@desargues.univ-lyon1.fr}

\Date{Received October 2, 2002; Accepted November 9, 2002}

\begin{abstract}
\noindent
We study non-trivial deformations of the natural action of the Lie
algebra $\mathrm{Vect}({\mathbb R}^n)$ on the space of differential forms on
${\mathbb R}^n$. We calculate abstractions for integrability of
infinitesimal multi-parameter deformations and determine the
commutative associative algebra corresponding to the miniversal
deformation in the sense of~\cite{ff}.
\end{abstract}

\section{Introduction}

The Lie algebra $\mathrm{Vect}(M)$ of vector fields on a smooth manifold
$M$ naturally acts on the space $\Omega(M)$ of differential
forms on $M$. Our goal is to study deformations of this action.

In this paper we restrict ourselves to the case $M={\mathbb R}^n$ with
$n\geq2$, and consider only those deformations which are given by
differentiable maps. We will use the framework of
Fialowski--Fuchs~\cite{ff} (see also~\cite{aalo}) and consider
(multi-parameter)
deformations over commutative associative algebras. We will
construct the miniversal deformation of Lie derivative and define
the commutative associative algebra related to this deformation.

The first step of any approach to the deformation theory consists
in the study of infinitesimal deformations. According to
Nijenhuis--Richardson~\cite{nr}, infinitesimal deformations are
classified by the first cohomology space
$H^1(\mathrm{Vect}(M);\mathrm{End}(\Omega(M)))$, while the obstructions for
integrability of infinitesimal deformations belong to the second
cohomo\-logy space $H^2(\mathrm{Vect}(M);\mathrm{End}(\Omega(M)))$. The first space
was calculated in~\cite{aal}; our task, therefore, is to calculate
the obstructions.

We will prove that all the multi-parameter deformations of the Lie
derivative are, in fact, infinitesimal, that is of first-order in
the parameters of deformation. The parameters have to satisfy a
system of quadratic relations; this defines interesting
commutative algebra. We will give concrete explicit examples of
the deformed Lie derivatives.

We mention that a similar problem was considered in~\cite{aalo} for the case of
symmetric contravariant tensor fields instead of differential forms.
In~\cite{ro1,ro2} a slightly different approach was used for
deformations of embeddings of Lie algebras.

\section{The main definitions}

In this section we introduce the algebra ${\mathcal D}(\Omega({\mathbb R}^n))$ of linear
differential operators on the space of differential forms on~${\mathbb R}^n$.
The Lie algebra $\mathrm{Vect}({\mathbb R}^n)$ has a natural embedding into
${\mathcal D}(\Omega({\mathbb R}^n))$ given by the Lie derivative.
We then describe the first cohomology space of $\mathrm{Vect}({\mathbb R}^n)$ with
coefficients in ${\mathcal D}(\Omega({\mathbb R}^n))$.

\subsection{Differential operators acting on differential forms}

The space of differential forms has a natural
(i.e.\ $\mathrm{Vect}({\mathbb R}^n)$-invariant) grading
\[
\Omega({\mathbb R}^n)=\bigoplus_{k=0}^n\Omega^k({\mathbb R}^n),
\]
where $\Omega^k({\mathbb R}^n)$ is the space of differential $k$-forms.
Let $\Omega_0^{n-k}({\mathbb R}^n)$ be the space of compactly supported
$(n-k)$-forms, there is a pairing
\begin{equation}
\label{Pair}
\Omega^k({\mathbb R}^n)\otimes\Omega_0^{n-k}({\mathbb R}^n)\to{\mathbb R}
\end{equation}
so that
$\Omega_0^{n-k}({\mathbb R}^n)
\subset
\left(
\Omega^k({\mathbb R}^n)
\right)^*$.

The space of differential operators on $\Omega({\mathbb R}^n)$ is an
associative algebra. As a $\mathrm{Vect}({\mathbb R}^n)$-module,
it is split into a direct sum of subspaces
\[
{\mathcal D}^{k,\ell}=
{\mathcal D}(\Omega^k({\mathbb R}^n),\Omega^\ell({\mathbb R}^n))
\]
with $k\leq\ell\leq n$. Obviously, one has a natural embedding
${\mathcal D}^{k,\ell}\subset
\Omega^\ell({\mathbb R}^n)\otimes
\left(
\Omega^k({\mathbb R}^n)
\right)^*$.

\subsection{The Lie derivative}

We introduce the notations that will be useful for
computations. We denote~$\xi^i$ the covector $dx^i$. A
differential $k$-form can be written
\[
\omega=\omega_{i_1\ldots i_k}(x)\,\xi^{i_1}\cdots\xi^{i_k}
\]
and $\xi^i$ are odd variables, i.e., $\xi^i\xi^j=-\xi^j\xi^i$.
The Lie derivative $L_X$ along a vector field~$X$
is given by the
differential operator on~$\Omega({\mathbb R}^n)$
\begin{equation}
\label{LieDer}
L_X=X^i\partial_{x^i}+
\frac{\partial{}X^i}{\partial x^j}\xi^j\partial_{\xi^i},
\end{equation}
where $\partial_{x^i}=\partial/\partial x^i$ and
$\partial_{\xi^i}=\partial/\partial \xi^i$; here and below
summation over repeated indices is understood.

\subsection{The first cohomology space}\label{FirstSect}

In the case $n\geq2$, the first cohomology space of
$\mathrm{Vect}({\mathbb R}^n)$ with coefficients in~${\mathcal D}^{k,\ell}$ was calculated
in~\cite{aal}. The result is as follows
\begin{equation}
\label{CohSpace}
H^1(\mathrm{Vect}({\mathbb R}^n);{\mathcal D}^{k,\ell})=
\left\{
\begin{array}{ll}
{\mathbb R},& \ell=k,\vspace{1mm}\\
{\mathbb R}^2,& \ell=k+1,\vspace{1mm}\\
{\mathbb R},& \ell=k+2,\vspace{1mm}\\
0,&\mbox{otherwise}.
\end{array}
\right.
\end{equation}

In the case $\ell=k$, this space is generated by the
cohomology class of the 1-cocycle
\begin{equation}
\label{eq10}
C_0(X)=\mathrm{Div}(X),
\end{equation}
where $X$ is a vector field and $\mathrm{Div}(X)$ is the divergence of $X$ with respect to
the standard volume form on ${\mathbb R}^n$.

In the case $\ell=k+1$, there are two non-trivial cohomology classes
corresponding to the 1-cocycles
\begin{equation}
\label{eq11}
C_{1}(X)=d\circ\mathrm{Div}(X),\qquad
\widetilde{C}_{1}(X)=\mathrm{Div}(X)\circ d,
\end{equation}
where $d$ is the de Rham differential.

Finally, in the case $\ell=k+2$, the cohomology space is spanned by
the class of the 1-cocycle
\begin{equation}
\label{eq12}
C_2(X)=d\circ\mathrm{Div}(X)\circ{}d.
\end{equation}

\begin{definition}
We denote $C_i^k\in{\mathcal D}^{k,k+i}$ the cocycle $C_i$
(and $\widetilde{C}_1$) restricted to~$\Omega^k({\mathbb R}^n)$:
\begin{equation}
\label{CocGrad} C^k_i(X):C_i(X)\mid_{\Omega^k({\mathbb R}^n)}, \qquad
i=0,1,2.
\end{equation}
Note that these elements
of $Z^1(\mathrm{Vect}({\mathbb R}^n);\Omega({\mathbb R}^n))$ are linearly independent
for every $k=0,\ldots,n-i$.
\end{definition}

\section{The general framework}

In this section we define deformations of Lie algebra
homomorphisms and introduce the notion of miniversal deformations
over commutative algebras.

\subsection{Infinitesimal deformations}

Let $\rho_0:{\mathfrak g}\to\mathrm{End}(V)$ be an action of a Lie algebra ${\mathfrak g}$
on a vector space $V$. When studying deformations of the
${\mathfrak g}$-action $\rho$, one usually starts with infinitesimal
deformations:
\[
\rho=\rho_0+tC,
\]
where $C:{\mathfrak g}\to\mathrm{End}(V)$ is a linear map and $t$ is a formal parameter.
The homomorphism condition
\[
[\rho(x),\rho(y)]=\rho([x,y]),
\]
where $x,y\in{\mathfrak g}$, is satisfied in order 1 in $t$ if and
only if $C$ is a 1-cocycle. The first cohomology  space
$H^1({\mathfrak g};\mathrm{End}(V))$ classifies infinitesimal deformations up to
equivalence (see, e.g., \cite{Fuc, nr}).

In the case when this space of cohomology is multi-dimensional, it
is natural to consider multi-parameter deformations. If
$\dim{H^1({\mathfrak g};\mathrm{End}(V))}=m$, then choose 1-cocycles
$C_1,\ldots,C_m$ representing a basis of $H^1({\mathfrak g};\mathrm{End}(V))$ and
consider the infinitesimal deformation
\begin{equation}
\label{InDefGen}
\rho=\rho_0+\sum_{i=1}^m{}t_i C_i,
\end{equation}
with independent parameters $t_1,\ldots,t_m$.

\newpage

In our study, the first cohomology space is determined by
the formula~(\ref{CohSpace}).
An infinitesimal deformation of the $\mathrm{Vect}({\mathbb R}^n)$-action
on $\Omega({\mathbb R}^n)$ is then of the form
${\mathcal L}_X=L_X+{\mathcal L}^{(1)}_X$, where $L_X$ is the Lie derivative
of differential forms along
the vector field $X$, and
\begin{equation}
\label{InfinDef} {\mathcal L}_X^{(1)}=\sum_{k=0}^n{}t_0^kC^k_0(X)+
\sum_{k=0}^{n-1}\left(t_1^kC^k_1(X)+
\widetilde{t}_1^k\widetilde{C}^k_1(X)\right)+
\sum_{k=0}^{n-2}{}t_2^kC^k_2(X),
\end{equation}
and where $t_0^k$, $t_1^k$, $\widetilde{t}_1^k$ and $t_2^k$ are
$4n$ independent parameters and the cocycles $C_0^k$, $C_1^k$,
$\widetilde{C}_1^k$ and $C_2^k$ are defined by formulae
(\ref{eq10})--(\ref{eq12}).

\subsection{Integrability conditions and deformations\\ over
commutative algebras}

Consider the problem of integrability of infinitesimal
deformations. Starting with the infinitesimal
deformation~(\ref{InDefGen}), one looks for a formal series
\begin{equation}
\label{BigDef}
\rho=
\rho_0+\sum_{i=1}^m{}t_iC_i+
\sum_{i,j}{}t_it_j\rho^{(2)}_{ij}+\cdots,
\end{equation}
where the highest-order terms $\rho^{(2)}_{ij}$,
$\rho^{(3)}_{ijk}$, $\ldots$
are linear maps from ${\mathfrak g}$ to $\mathrm{End}(V)$ such that
\[
\rho:{\mathfrak g}\to\mathrm{End}(V)\otimes{\mathbb C}[[t_1,\ldots,t_m]]
\]
satisfy the homomorphism condition in any order in
$t_1,\ldots,t_m$.

However, quite often the above problem has no solution.
Following~\cite{ff} and~\cite{aalo}, we will impose extra algebraic
relations on the parameters $t_1,\ldots,t_m$. Let ${\mathcal R}$ be an
ideal in ${\mathbb C}[[t_1,\ldots,t_m]]$ generated by some set of
relations, the quotient
\begin{equation}
\label{TrivAlg}
 {\mathcal A}={\mathbb C}[[t_1,\ldots,t_m]]/{\mathcal R}
\end{equation}
is a commutative associative algebra with unity, and one can speak
of a deformations with base~${\mathcal A}$, see \cite{ff} for details. The
map (\ref{BigDef}) sends ${\mathfrak g}$ to $\mathrm{End}(V)\otimes{\mathcal A}$.

\begin{example}
\label{Example}
Consider the ideal ${\mathcal R}$ generated by all the quadratic monomials
$t_it_j$. In this case
\begin{equation}
\label{InfAlg} {\mathcal A}={\mathbb C}\oplus{\mathbb C}^m
\end{equation}
and any deformation is of the form (\ref{InDefGen}). In this case
any infinitesimal deformation becomes a deformation with base
${\mathcal A}$ since $t_it_j=0$ in ${\mathcal A}$, for all $i,j=1,\ldots,m$.
\end{example}

Given an infinitesimal deformation (\ref{InDefGen}), one can
always consider it as a deformation with base~(\ref{InfAlg}). Our
aim is to find ${\mathcal A}$ which is big as possible, or, equivalently,
we look for relations on $t_1,\ldots,t_m$ which are necessary and
sufficient for integrability (cf.~\cite{aalo}).

\subsection{Equivalence and the miniversal deformation}

The notion of equivalence of deformations over commutative
associative algebras has been considered in~\cite{ff}.

\setcounter{definition}{1}
\begin{definition}
Two deformations $\rho$ and $\rho'$ with the same base ${\mathcal A}$
are called equivalent if there exists an inner automorphism
$\Psi$ of the associative algebra $\mathrm{End}(V)\otimes{\mathcal A}$ such that
\[
\Psi\circ\rho=\rho'
\]
and such that
$\Psi({\mathbb I})={\mathbb I}$,
where ${\mathbb I}$ is the unity of the algebra $\mathrm{End}(V)\otimes{\mathcal A}$.
\end{definition}

The following notion of miniversal deformation is fundamental.
It assigns to a~${\mathfrak g}$-modu\-le~$V$
a canonical commutative associative algebra ${\mathcal A}$
and a canonical deformation with base~${\mathcal A}$.

\begin{definition}
A deformation~(\ref{BigDef}) with base ${\mathcal A}$
is called miniversal, if for any other deformation $\rho'$
with base ${\mathcal A}'$ there exists a unique homomorphism
$\psi:{\mathcal A}\to{\mathcal A}'$ satisfying $\psi(1)=1$,
such that
\[
\rho'=(\mathrm{Id}\otimes\psi)\circ\rho.
\]
This definition
does not depend on the choice of the basis $C_1,\ldots,C_m$.
\end{definition}

The miniversal deformation corresponds to the smallest ideal
${\mathcal R}$. We refer to~\cite{ff} for a~construction of miniversal
deformations of Lie algebras and to~\cite{aalo} for miniversal
deformations of ${\mathfrak g}$-modules.

\section{The main result}

In this section we obtain the integrability conditions for the
infinitesimal deformation~(\ref{InfinDef}). The main result of
this paper is the following
\begin{theorem}
\label{Main}
(i)
The following $4n-4$ conditions
\begin{align}
& R^k_1(t)=t^k_0\widetilde{t}^k_1+
t_0^{k+1}t_1^k=0, &&
\mbox{for}\quad  k=0,\ldots,n-1,& \nonumber\\
 & R^k_2(t)=t_0^kt_2^k+t_1^{k+1}t_1^k=0,
 &&  \mbox{for}\quad k=0,\ldots,n-2,& \nonumber\\
& \widetilde{R}^k_2(t)=
t_0^{k+2}t_2^k+\widetilde{t}^{k+1}_1\widetilde{t}^k_1=0,
  && \mbox{for}\quad k=0,\ldots,n-2,& \nonumber\\
& R^k_3(t)=t_1^{k+2}t_2^k+\widetilde{t}^k_1t_2^{k+1}=0,
&& \mbox{for}\quad k=0,\ldots,n-3&
\label{CondInt}
\end{align}
are necessary and sufficient for integrability of the
infinitesimal deformation~(\ref{InfinDef}).

(ii) Any
formal deformation of the $\mathrm{Vect}({\mathbb R}^n)$-action on $\Omega({\mathbb R}^n)$
is equivalent to a deformation of order 1, that is,
to a deformation given by~(\ref{InfinDef}).
\end{theorem}

The commutative algebra defined by relations~(\ref{CondInt})
corresponds to the miniversal deformation of the Lie derivative
$L_X$. Note that the commutative algebra
defined in Theorem~\ref{Main} is infinite-dimensional.

We start the proof of Theorem~\ref{Main}. It consists in two
steps. First, we compute explicitly the obstructions for existence
of the second-order term, this will prove that relations~(\ref{CondInt})
are necessary. Second we show that under
relations~(\ref{CondInt}) the highest-order terms of the deformation can be
chosen identically zero, so that relations~(\ref{CondInt}) are
indeed sufficient.

\subsection{Computing the second-order Maurer--Cartan equation}

Assume that the infinitesimal deformation~(\ref{InfinDef}) can be
integrated to a formal deformation
\[
{\mathcal L}_X=L_X+{\mathcal L}^{(1)}_X+{\mathcal L}^{(2)}_X+\cdots,
\]
where ${\mathcal L}^{(1)}_X$ is given by (\ref{InfinDef}) and ${\mathcal L}^{(2)}_X$
is a quadratic polynomial in $t$ with coefficients in~${\mathcal D}(\Omega({\mathbb R}^n))$.
We compute the conditions for the second-order terms~${\mathcal L}^{(2)}$.
Consider the quadratic terms of the homomorphism condition
\begin{equation}
\label{HomomCond}
[{\mathcal L}_X,{\mathcal L}_Y]={\mathcal L}_{[X,Y]}.
\end{equation}
The left hand side of (\ref{HomomCond}) contains the operators
from $\Omega^k({\mathbb R}^n)$ to $\Omega^{k+i}({\mathbb R}^n)$, for every
$k=0,\ldots,n$ and $i=1,2,3$. Condition~(\ref{HomomCond}) must be
satisfied in any order of~$i$ independently.

Collecting the terms with $i=1$, one readily gets:
\[
R^k_1(t)\gamma_1(X,Y),
\]
where $\gamma_1$
is the bilinear skew-symmetric map from $\mathrm{Vect}({\mathbb R}^n)$ to
${\mathcal D}^{k,k+1}$ defined by
\begin{equation}
\label{2-Coc1}
\gamma_1(X,Y)=\left(
\mathrm{Div}(X)d(\mathrm{Div}(Y))-\mathrm{Div}(Y)d(\mathrm{Div}(X))
\right)\wedge.
\end{equation}
The terms with $i=2$ in (\ref{HomomCond}) are as follows
\[
R^k_2(t)\gamma_2(X,Y)+\widetilde{R}^k_2(t)\,\widetilde{\gamma}_2(X,Y),
\]
where $\gamma_2$ and $\widetilde{\gamma}_2$ are the bilinear
skew-symmetric maps from $\mathrm{Vect}({\mathbb R}^n)$ to
${\mathcal D}^{k,k+2}$ defined by
\begin{equation}
\label{2-Coc2}
\gamma_2(X,Y)=\left(
\mathrm{Div}(X)d(\mathrm{Div}(Y))-\mathrm{Div}(Y)d(\mathrm{Div}(X))
\right)\wedge{}d
\end{equation}
and
\begin{equation}
\label{2-Coc2tilde}
\widetilde{\gamma}_2(X,Y)=\left(
d(\mathrm{Div}(X))\wedge{}d(\mathrm{Div}(Y))-d(\mathrm{Div}(Y))\wedge{}d(\mathrm{Div}(X))
\right)\wedge.
\end{equation}
Finally, the terms with $i=3$ in (\ref{HomomCond}) are as follows
\[
R^k_3(t)\gamma_3(X,Y),
\]
where $\gamma_3$ is the bilinear
skew-symmetric map from $\mathrm{Vect}({\mathbb R}^n)$ to
${\mathcal D}^{k,k+3}$ defined by
\begin{equation}
\label{2-Coc3}
\gamma_3(X,Y)=\left(
d(\mathrm{Div}(X))\wedge{}d(\mathrm{Div}(Y))-
d(\mathrm{Div}(Y))\wedge{}d(\mathrm{Div}(X))
\right)\wedge{}d,
\end{equation}
where the coefficients $R^k_1(t)$,
$R^k_2(t)$, $\widetilde{R}^k_2$ and
$R^k_3$ are precisely the quadratic polynomial from~(\ref{CondInt}).

The homomorphism condition~(\ref{HomomCond}) gives for the second-order
terms the following (Maurer--Cartan) equation
\[
\delta({\mathcal L}^{(2)})=-\frac 12 [\![{\mathcal L}^{(1)},{\mathcal L}^{(1)}]\!],
\]
where $\delta$ is the Chevalley--Eilenberg differential and
$[\![\,,\,]\!]$ stands for the cup-product of 1-cocycles, so that
the right hand side of this expression is automatically a
2-cocycle. In our case, we obtain explicitly:
\begin{equation}
\label{MC2}
\delta({\mathcal L}^{(2)})=R^k_1(t)\gamma_1+
R^k_2(t)\gamma_2+\widetilde{R}^k_2(t)\,\widetilde{\gamma}_2+
R^k_3(t)\gamma_3.
\end{equation}
The bilinear maps $\gamma_1$,
$\gamma_2$, $\widetilde{\gamma}_2$ and $\gamma_3$ are
2-cocycles on $\mathrm{Vect}({\mathbb R}^n)$. For existence of solutions of the
equation~(\ref{MC2}) it is necessary and sufficient that the right
hand side be a coboundary.

\subsection{Obstructions for integrability}

We denote
$\gamma^k_1$, $\gamma^k_2$, $\widetilde{\gamma}^k_2$ and $\gamma^k_3$ the
2-cocycles on $\mathrm{Vect}({\mathbb R}^n)$ such that
\[
\gamma^k_i(X,Y)=\gamma_i(X,Y)\big|_{\Omega^k({\mathbb R}^n)}
\]
for every $X,Y\in\mathrm{Vect}({\mathbb R}^n)$.
In order to solve the equation~(\ref{MC2}), we now have to study
the cohomology classes of these 2-cocycles in
$H^2(\mathrm{Vect}({\mathbb R}^n);{\mathcal D}^{k,k+i})$.

\setcounter{proposition}{1}
\begin{proposition}
\label{ProNontriv}
Each of the 2-cocycles:
\begin{gather*}
\gamma^k_1,\phantom{\ \ \widetilde{\gamma}^k_2,} \qquad \mbox{for}\quad k=0,\ldots,n-1,\\
\gamma^k_2,\ \ \widetilde{\gamma}^k_2, \qquad \mbox{for}\quad k=0,\ldots,n-2
\quad\mbox{and}\quad n\geq2,  \\
\gamma^k_3,\phantom{\ \ \widetilde{\gamma}^k_2,} \qquad \mbox{for}\quad k=0,\ldots,n-3
\quad\mbox{and}\quad n\geq3
\end{gather*}
define non-trivial cohomology class. Moreover, these cohomology
classes are linearly independent.
\end{proposition}

\begin{proof}
Assume that, for some differential 1-cohain $b$ on $\mathrm{Vect}({\mathbb R}^n)$
with coefficients in ${\mathcal D}^{k,k+i}$, one has $\gamma^k_i=\delta(b)$.
The general form of such a cochain is
\[
b(X)=
b^{i_1\ldots{}i_sj_1\ldots{}j_tm_1\ldots{}m_r}_{\ell\,a_1\ldots{}a_r}(x)\,
\frac{\partial^s\,X^\ell}{\partial x^{i_1}\cdots\partial x^{i_s}}\,
\xi^{a_1}\cdots\xi^{a_r}
\partial_{x^{j_1}}\cdots\partial_{x^{j_t}}
\partial_{\xi^{m_1}}\cdots\partial_{\xi^{m_u}}.
\]
We will need the following:

\setcounter{lemma}{2}
\begin{lemma}
The condition $\delta(b)=\gamma^k_i$ implies that the coefficients
$b^{i_1\ldots{}i_sj_1\ldots{}j_tm_1\ldots{}m_r}_{\ell\,a_1\ldots{}a_r}$
are constants.
\end{lemma}

\begin{proof}
The condition $\gamma^k_i=\delta(b)$ reads
\begin{equation}
\label{delta}
\gamma^k_i(X,Y)=[L_X,b(Y)]-[L_Y,b(X)]-b([X,Y]).
\end{equation}
We choose a constant vector field $Y=\partial_{x^i}$ and prove that
$L_Y(b)=0$. Indeed, one has
\[
(L_Y(b))(X)=[L_Y,b(X)]-b([Y,X]).
\]
Since $s\geq2$ in the expression of $b$, it follows that $b(Y)=0$,
and thus the last equality gives
\[
(L_Y(b))(X)=[L_Y,b(X)]-b([Y,X])-[L_X,b(Y)]=(\delta b)(Y,X).
\]
By assumption, $\delta b=\gamma^k_i$, from the explicit formula for
$\gamma^k_i$, see (\ref{2-Coc1})--(\ref{2-Coc3}) one obtains
$\gamma^k_i(Y,X)=0$ for all $X$, since $\mathrm{Div}(Y)=0$. Therefore,
$L_Y(b)=0$.

Lemma~4.3 is proved.
\end{proof}

It is now easy to check that the equation (\ref{delta}) has no
solutions: using the formula~(\ref{LieDer}), we see that the
expressions $\mathrm{Div}(X)$ and $\mathrm{Div}(Y)$ never appear in the right hand
side of~(\ref{delta}). This is a contradiction with the assumption
$\gamma^k_i=\delta(b)$.

Proposition~\ref{ProNontriv} is proved.\end{proof}

Proposition~\ref{ProNontriv} implies that the equation~(\ref{MC2})
has solutions if and only if the quadratic polynomials
$R^k_1(t)$, $R^k_2(t)$,
$\widetilde{R}^k_2(t)$ and $R^k_3(t)$
vanish simultaneously. We thus proved that the
conditions~(\ref{CondInt}) are, indeed, necessary.

\subsection{Integrability conditions are sufficient}

To prove that the conditions (\ref{CondInt}) are sufficient, we
will find explicitly a deformation of~$L_X$, whenever the
conditions~(\ref{CondInt}) are satisfied. The solution ${\mathcal L}^{(2)}$
of~(\ref{MC2}) can be chosen identically zero. Choosing the
highest-order terms ${\mathcal L}^{(m)}$ with $m\geq3$, also identically
zero, one obviously obtains a deformation (which is of order 1 in~$t$).

Theorem~\ref{Main}, part (i) is proved.

The solution ${\mathcal L}^{(2)}$ of (\ref{MC2}) is defined up to a
1-cocycle and it has been shown in~\cite{ff,aalo} that different
choices of solutions of the Maurer--Cartan equation correspond to
equivalent deformations. Thus, one can always reduce ${\mathcal L}^{(2)}$
to zero by equivalence. Then, by recurrence, the highest-order
terms ${\mathcal L}^{(m)}$ satisfy the equation $\delta({\mathcal L}^{(m)})=0$ and can
also be reduced to the identically zero map. This completes the
proof of part~(ii).

Theorem \ref{Main} is proved.

\section{An open problem}

It seems to be an interesting open problem to compute the full
cohomology ring \linebreak $H^*(\mathrm{Vect}({\mathbb R}^n); {\mathcal D}^{k,\ell})$. The only
complete result here concerns the first cohomology space,
see~\cite{aal}. Proposition~\ref{ProNontriv} provides a lower bound
for the dimension of the second cohomology space. We formulate

\begin{conjecture}
The space of second cohomology of $\mathrm{Vect}({\mathbb R}^n)$
with coefficients in the space of differential operators
on $\Omega({\mathbb R}^n)$ has the following structure
\[
H^2(\mathrm{Vect}({\mathbb R}^n);{\mathcal D}^{k,\ell})=
\left\{
\begin{array}{ll}
{\mathbb R},& \ell=k+1,\vspace{1mm}\\
{\mathbb R}^2,& \ell=k+2,\vspace{1mm}\\
{\mathbb R},& \ell=k+3,\vspace{1mm}\\
0,&\mbox{otherwise}
\end{array}
\right.
\]
and spanned by the $2$-cocycles $\gamma^k_1$,
$\gamma^k_2$, $\tilde{\gamma}^k_2$
and $\gamma^k_3$ given by (\ref{2-Coc1})--(\ref{2-Coc3}).
\end{conjecture}

\section{A few examples of deformations}

We give some examples of deformations of the
$\mathrm{Vect}({\mathbb R}^n)$-action (\ref{LieDer}) on $\Omega({\mathbb R}^n)$ which
are simpler than the general case described by Theorem~\ref{Main}.

\begin{example}
The first example is a 1-parameter deformation
\[
{\mathcal L}_X=L_X+tC_0(X),
\]
that is, we put $t_0^k=t$ and $t_1^k=\widetilde{t}_1^k=t_2^k=0$.
This deformation has a geometric meaning. More precisely,
we consider the $\mathrm{Vect}({\mathbb R}^n)$-module of (generalized) tensor fields
$\Omega({\mathbb R}^n)\otimes_{C^\infty({\mathbb R}^n)}
{\mathcal F}_t({\mathbb R}^n)$,
where ${\mathcal F}_t({\mathbb R}^n)$ is the space of tensor densities of degree~$t$
on ${\mathbb R}^n$.
\end{example}

\begin{example}
Another simple example of a 1-parameter deformation is given by
\[
{\mathcal L}_X=L_X+tC_2(X).
\]
The geometric meaning of this deformation is not clear.
\end{example}

\setcounter{proposition}{2}
\begin{proposition}
The above two examples are the only 1-parameter
deformations of the form
\[
{\mathcal L}_X=L_X+t\left(\alpha_0 C_0(X)+\alpha_1C_1(X)
+\widetilde{\alpha}_1\widetilde{C}_1(X)
+\alpha_2C_2(X)\right),
\]
that is, with $t_0^k=\alpha_0 t$,
$t_1^k=\alpha_1t$,
$\widetilde{t}_1^k=\widetilde{\alpha}_1t$ and $t_2^k=\alpha_2t$ for all $k$.
\end{proposition}

\begin{proof}
Straightforward from the relations~(\ref{CondInt}).
\end{proof}

\setcounter{example}{3}
\begin{example}
We give an example of 3-parameter deformations in the
two-dimen\-sio\-nal case, i.e., for $n=2$
\[
{\mathcal L}_X=L_X+
t_0^2\left(C^2_0(X)-\widetilde{C}_1^0(X)\right)+
t_1^0 C_1^0(X)
+\widetilde{t}_1^1\left(\widetilde{C}_1^1(X)+C^0_2(X)\right),
\]
where $t_0^2$, $t_1^0$ and $\widetilde{t}_1^1$ are independent parameters.
\end{example}

One can construct a great number of examples of deformations with
3 independent parameters in the multi-dimensional case; it would be
interesting to understand their geometric meaning.

\subsection*{Acknowledgments}

It is our pleasure to acknowledge enlightening
discussions with C~Duval, P~Lecomte, D~Leites and C~Roger.

\label{ovsienko-lastpage}

\end{document}